\documentclass[11pt, a4paper]{amsart}
\usepackage{amsmath,amssymb,amscd,amsfonts}

\newtheorem{thm}{Theorem}[section]
\newtheorem{lem}[thm]{Lemma}
\newtheorem{rem}[thm]{Remark}
\newtheorem{prop}[thm]{Proposition}
\newtheorem{cor}[thm]{Corollary}

\newtheorem{step}{Step}
\newcommand{\Def}{\text{\rm Def}}
\newcommand{\Proj}{\text{\rm Proj}}

\newcommand{\Aut}{\text{\rm Aut}}

\newcommand{\OO}{{\mathcal {O}}}

\newcommand{\inv}{^{-1}}
\newcommand{\al}{\alpha}
\newcommand{\be}{\beta}
\newcommand{\ga}{\gamma}
\newcommand{\Ga}{\Gamma}

\newcommand{\epsi}{\epsilon}
\newcommand{\fie}{\varphi}
\newcommand{\si}{\sigma}
\newcommand{\Si}{\Sigma}

\newcommand{\C}{{\mathbb C}}
\newcommand{\R}{{\mathbb R}}
\newcommand{\Q}{{\mathbb Q}}
\newcommand{\Z}{{\mathbb Z}}
\newcommand{\F}{{\mathbb F}}
\newcommand{\pp}{{\mathbb{P}}}

\newcommand{\CC}{{\mathcal{C}}}
\newcommand{\DD}{{\mathcal{D}}}
\newcommand{\FF}{{\mathcal{F}}}

\newcommand{\MM}{{\mathcal{M}}}

\newcommand{\XX}{{\mathcal{X}}}
\numberwithin{equation}{section}

\title[The classification of double planes\dots]{The classification of double
planes of general type with $K^2=8$ and
$p_g=0$}
\author{Rita Pardini}
\date{}
\begin{document}

\begin{abstract}We study minimal   {\em double planes} of general type with $K^2=8$ and
$p_g=0$, namely pairs $(S,\si)$,  where $S$ is a minimal complex algebraic surface of general
type  with
$K^2=8$ and $p_g=0$ and $\si$ is an automorphism of $S$ of order 2 such that the quotient
$S/\si$ is a rational surface. We prove that $S$ is a free quotient $(F\times C)/G$, where
$C$ is a curve, $F$ is an hyperelliptic curve, $G$ is a finite group that acts faithfully on
$F$ and $C$, and $\si$ is induced by the automorphism $\tau\times Id$ of $F\times C$, $\tau$
being the hyperelliptic involution of $F$. We describe all
the $F$, $C$  and $G$ that occur: in this way we obtain 5
families of surfaces with $p_g=0$ and $K^2=8$, of which we
believe only one was previously known.  

Using our classification we are able to
give an alternative  description of these surfaces as double covers of the plane, thus
recovering a construction proposed by Du Val. In addition we  study the geometry of the
   subset of the moduli space of surfaces of general type  with $p_g=0$ and $K^2=8$ that
admit a double plane structure.
\smallskip

\noindent 2000 Mathematics Classification: 14J29.
\end{abstract}
\maketitle
\section{Introduction}

The last section of the paper \cite{nodi}  is concerned with the study  of pairs
$(S,\si)$,  where $S$ is a minimal complex projective  surface of general type
with $p_g=0$ and $K^2=8$  and
$\si$ is an automorphism of $S$ of order 2. Here we give a complete description in the case
that
$(S,
\si)$ is a {\em double plane}, i.e. that the quotient surface $S/\si$ is rational.

First of all, using the results of \cite{nodi}, we prove that $S$ is a free
quotient of a product of curves:
\begin{thm}\label{A} 
Let $S$ be a minimal complex projective surface of general type   with $K^2_S=8$ and
$p_g(S)=0$.

There exists an automorphism $\si$ of $S$ of order $2$ such that $S/\si$ is rational if and
only if there exist a curve $C$, an hyperelliptic curve $F$ of genus 3 or 5  and a finite
group
$G$ such that:
\begin{itemize}
\item[a)] $G$ acts faithfully on $F$ and $C$  and the diagonal action of $G$
on $F\times C$ is free;
\item[b)] $|G|=(g(F)-1)(g(C)-1)$;
\item[c)] $C/G$ and $F/G$ are rational curves;
\item[d)] $S=(F\times C)/G$  and $\si$ is the involution induced by
$\tau\times Id$, where $\tau$ denotes the hyperelliptic involution of $F$. 
\end{itemize}
\end{thm}

Theorem \ref{A} reduces the classification of double planes of general type  with $p_g=0$
and $K^2=8$   to the classification of the triples
$F$,
$C$,
$G$  that satisfy properties  a)--c) of its statement. 
The latter problem can be attacked by remarking that the $G-$action on the hyperelliptic
curve $F$ descends to a (possibly not  faithful) $G-$action on $F/\tau=\pp^1$ that preserves
the image of the set of  Weierstrass points of $F$. We  combine this observation with the
well known classification of the finite group actions on the projective line and, exploiting
also  some other geometrical constraints of the problem, we obtain a complete classification.
There are 5 types of minimal  double planes with $p_g=0$ and $K^2=8$, whose description takes
up  all of section
\ref{sezexamples}. We list here only the numerical invariants for each type:
\begin{thm}\label{B} Let $(S, \si)$ be a minimal double plane with $p_g=0$ and $K^2=8$. Then,
referring to the notation of Theorem \ref{A}, $(S,\si)$ belongs to one and only one of the
following types:
\begin{itemize}
\item[Ia:] $g(F)=3$, $g(C)=5$, $G=\Z_2^3$;
\item[Ib:] $g(F)=3$, $g(C)=9$, $G=\Z_2\times D_4$;
\item[Ic:] $g(F)=3$, $g(C)=13$, $G=S_4$;
\item[Id:] $g(F)=3$, $g(C)=25$, $G=\Z_2\times S_4$;
\item[II:] $g(F)=5$, $g(C)=16$, $G=A_5$.
\end{itemize}
\end{thm}

While double planes of type Ia already appear in \cite{mp}, as far as we know the
remaining 4 types give 4 new families of surfaces of general type with $p_g=0$
and $K^2=8$. It may be  worth remarking that all known surfaces with these
invariants are free quotients of product of curves (a construction  suggested
first by Beauville). It would be interesting to know whether these are actually
the only ones.

Another byproduct of Theorem \ref{B}, which was partly the
motivation of this work, is the classification of the
minimal  surfaces of general type with $p_g=0$, $K^2=8$ and
non birational bicanonical map.
\begin{cor}
Let $S$ be a minimal projective complex surface of general type with $p_g=0$ and $K^2=8$.

The bicanonical map of $S$ is non birational if and only if there is an automorphism $\si$
of $S$ such that $(S,\si)$ is a double plane of type Ia, Ib, Ic or Id (cf. Thm.
\ref{B}).
\end{cor}

 Furthermore,  using Theorem \ref{B}
we are able to give an alternative description of these surfaces  as minimal resolution of
normal double covers of the projective plane $\pp^2$ (``plane models''), thus recovering a
construction suggested by Du Val (cf.
\S
\ref{sezmodels}). We also study the subset
 of the moduli space of surfaces with $p_g=0$ and $K^2=8$ consisting of surfaces  that admit
a double plane construction (cf. \S \ref{sezmoduli}). In particular we show that surfaces
with non birational bicanonical map are a closed and open subset of the moduli
space  which is the union of four irreducible connected components.

The paper goes as follows: \S 2 contains the proof of Theorem \ref{A};  \S 3
contains the description of the 5 types of double planes; \S 4 is the technical heart of the
paper and contains the proof that the double planes with $K^2=8$ and $p_g=0$ are precisely
those given in \S 3; \S 5 contains the description of the plane models;
finally \S 6 is devoted to the study of the moduli.
\bigskip 

\noindent{\bf Acknowledgments.} The starting point of this work  is  the last section of
the joint paper \cite{nodi} and I am indebted   to both
coauthors, Igor  Dolgachev and Margarida Mendes Lopes, for many interesting discussions.
In particular, the first section should be regarded as joint work with
Margarida Mendes Lopes. Finally, I am  grateful to Barbara Fantechi for her
advice on the deformation theory of the last section.
\bigskip

\noindent{\bf Notation and conventions.} We work over the field of complex numbers. By {\em
surface\ } we mean  a   projective surface and by {\em curve} we mean a smooth
projective  curve.  A {\em map\
} is a  rational map and a {\em morphism\ } is a regular map. If a group $G$ acts on a set
$X$ we denote by
$X^G$ the subset of elements fixed by $G$. 

 We use standard notation of algebraic geometry; we
just recall here the notation for the invariants of a surface $S$: $K_S$
is the {\em canonical class}, $p_g(S)=h^0(S,K_S)$ the {\em geometric
genus\ },  $q(S)=h^1(S,\OO_S)$ the {\em irregularity} and $\chi(S)=1+p_g(S)-q(S)$ the {\em
Euler characteristic}.

\section{Rational involutions}\label{sezinvolutions}
Let $S$ be a smooth projective surface.
An {\it involution} $\si$ of   $S$ is an automorphism of $S$ of order 2.  We say that a
 map $f\colon S\to Y$ is {\it composed with}  $\si$\ if
$f\circ\si=f$.
 
An involution
$\si$ is called {\em rational}  if the quotient surface
$\Si:=S/\si$ is rational. 
We define a {\em double
plane} as a pair
$(S,\si)$, where
$S$ is a smooth projective surface  and $\si$ is a rational involution of $S$. 
We say that a   double plane $(S,\si)$ is  minimal, of general
type, has $p_g=0$, \dots if the surface
$S$ is minimal, of general type, has $p_g=0$\dots 

Isomorphism of double planes is defined in the
obvious way. 
\smallskip

The following proposition explains how to construct  examples of minimal double planes
$(S,\si)$ by taking free quotients of product of curves.
\begin{prop}\label{costruzione}
Let $F$ be an hyperelliptic curve and let $\tau$ be the hyperelliptic involution
of $F$,  let
$C$ be a curve,   and let
$G$ be  a finite group that acts faithfully on $F$ and on $C$. Assume that:
\begin{itemize}
\item[a)]
$C/G$ is rational;
\item[b)] 
the diagonal action of
$G$ on
$F\times C$ is free.
\end{itemize} Set
$S:=(F\times C)/G$ and denote by $\si$ the involution of $S$ induced by
$\tau\times Id$.

Then $(S, \si)$ is a  minimal
double plane with invariants: 
$$ \chi(S)=\frac{(g(F)-1)(g(C)-1)}{|G|};\quad q(S)=g(F/G).$$
$$K_S^2=\frac{8(g(F)-1)(g(C)-1)}{|G|}.$$

Moreover, $S$ is a minimal surface of general type  iff $g(F),g(C)\ge 2$.
\end{prop}
\begin{proof}Consider the quotient map $\psi\colon F\times C\to S$. By 
assumption b), the map  $\psi$  is   \'etale   of degree equal to $|G|$,
hence
$S$ is smooth and  we have
$\chi(F\times C)=|G|\chi(S)$ and $K_{F\times C}=\psi^*K_S$.  The formulas
for $\chi(S)$, $K^2_S$ follow easily from this remark. In addition,
$S$ is a minimal surface of general type iff $K_S$ is nef and big, iff
$K_{F\times C}$ is nef and big, iff $C$ and $F$ have genus $\ge 2$. 
The irregularity of $S$ is equal to the dimension of the
$G-$invariant subspace of $H^0(\Omega^1_{F\times C})=H^0(\omega_F)\oplus
H^0(\omega_C)$. Since $G$ acts separately on $F$ and $C$, one has
$q(S)=g(C/G)+g(F/G)=g(F/G)$. We set, as usual, $\Si:=S/\si$. The second
projection $F\times C\to C$ induces a pencil
$\Si\to C/G=\pp^1$ whose general fibre is isomorphic to $\pp^1$, hence $\Si$
is rational and $(S,\si)$ is a double plane.
\end{proof}

 The possible types of involutions of a minimal surface of general type  $S$
with
$p_g=0$ and
$K^2=8$ are described in Theorem 4.4 of \cite{nodi}. In the case of double
planes this description can be made more precise, showing in
particular that all minimal  double planes with $p_g=0$ and $K^2=8$ arise as
in Proposition
\ref{costruzione}.

\begin{thm}\label{struttura}
Let $(S,\si)$ be a minimal double plane with $K^2_S=8$ and $p_g(S)=0$.

Then
there exist a curve $C$, an hyperelliptic curve $F$ with hyperelliptic
involution $\tau$ and a finite group
$G$ such that:
\begin{itemize}
\item[a)] $G$ acts faithfully on $F$ and $C$  and the diagonal action of $G$
on $F\times C$ is free;
\item[b)] $C/G$ and $F/G$ are rational curves;
\item[c)] $S=(F\times C)/G$  and $\si$ is the involution induced by
$\tau\times Id$. 
\end{itemize}

 \noindent  If we  denote by 
$k$ be the number of isolated fixed points of $\si$, then  the numerical
possibilities for $k$, $g(F)$ and $g(C)$ are the following:

\begin{itemize}
\item[i)]  $k=12$, $g(F)=3$,   $|G|=2(g(C)-1)$;
\item[ii)]  $k=10$, $g(F)=5$,   $|G|=4(g(C)-1)$.
\end{itemize}
\end{thm}
\begin{proof}By Theorem 4.4 of \cite{nodi}, there are the following
possibilities for $(S,\si)$:
\begin{itemize}
\item[i)] $k=12$,  $S$ has a  pencil $p\colon S\to\pp^1$ of hyperelliptic
curves of genus 3 with 6 double fibres and $\si$ restricts to the hyperelliptic
involution on the general fibre  $F$ of $p$;
\item[ii)]  $k=10$,  $S$ has a pencil $p\colon S\to\pp^1$ of hyperelliptic
curves of genus 5 with 5 double fibres and $\si$ restricts to the hyperelliptic
involution on the general fibre $F$ of $p$.
\end{itemize}
Consider case i) first. Let $h\colon B\to \pp^1$ be the double cover
branched on the 6 points corresponding to the double fibres of $p$. $B$ is a
smooth curve of genus 2.  By taking base change and normalization, one gets a
diagram:
\[\begin{CD}
X@>\tilde{h}>>S\\
@V\tilde{p}VV  @VVpV\\
B@> h >> \pp^1
\end{CD}\]
where $\tilde{h}\colon X\to S$ is an \'etale double cover and $\tilde{p}\colon 
X\to B$ is a pencil with general fibre isomorphic to $F$. One has
$16=K^2_X=8(g(F)-1)(g(B)-1)$, hence by \cite{appendice} the fibration  $\tilde{p}$ is 
smooth and isotrivial. The above diagram shows that $p$ is also
isotrivial, that the 6 double fibres have smooth support and they are the
only singular fibres of $p$. So, in the terminology of \cite{serrano}, $p\colon
S\to\pp^1$ is a {\em quasi-bundle} and by \S 1 of \cite{serrano} (cf. also \cite{serranoC})
there exist a curve
$C$ and a finite group $G$ that acts faithfully on $C$ and $F$ in such a way that the
diagonal action on $F\times C$ is free and 
$S$ is isomorphic to
$(F\times C)/G$. The pencil $p$ is
induced by the second projection $F\times C\to C$ and the involution $\si$ is
induced by $\tau\times Id$, where $\tau$ is the hyperelliptic involution of
$F$. Since $S$ is regular, we have $g(F/G)=g(C/G)=0$ (cf. 
 proof of Proposition \ref{costruzione}). The formula for $|G|$ follows
from Proposition \ref{costruzione}. 

 If
$S$ is as in case ii), then we let
$h\colon B\to
\pp^1$  be a
$\Z_2^2-$cover branched over the 5 points corresponding to the double fibres
of $p$ (the existence of such a cover can be easily shown by using Theorem 1.2
of \cite{ritaabel}). $B$ is again a smooth curve of genus 2 and  one argues
exactly as before.
\end{proof}

Minimal surfaces of general type with $p_g=0$, $K^2=8$ and non birational
bicanonical map are a special  instance of double planes. In fact they are
precisely the double planes of case i) of Theorem \ref{struttura}, as 
explained below.
\begin{cor}\label{bica}Let $S$ be a smooth minimal surface of general type with $p_g(S)=0$
and
$K^2_S=8$. 

 The bicanonical map
$\fie$ of
$S$ is not birational onto its image iff there exist a curve $C$, an hyperelliptic curve
$F$ of genus
$3$ and a group
$G$ of order
$2(g(C)-1)$ such that:
\begin{itemize}
\item[a)]  $G$ acts faithfully on $C$ and on $F$ and the diagonal action on
$F\times C$ is free;
\item[b)]
$S=(F\times C)/G$ 
\end{itemize}
In this case   $\fie$ has degree $2$  and it is   composed with the
involution
$\si$ of
$S$ 
 induced by
$\tau\times Id$, where $\tau$ is the bicanonical involution of $F$. The
involution  
$\si$ is rational and it   has 12 isolated fixed points.
\end{cor}
\begin{proof}
By Theorem 1.1 of \cite{mp} $\fie$ is not birational iff it has degree 2, hence
$\fie$ is not birational iff it is composed with an involution  $\si$ of $S$.
By the results of \cite{mp3} this happens iff $\si$ has $12$ isolated
fixed points. Finally, by Theorem 3 of \cite{xiaodeg} if $\fie$ has
degree 2 then
$\fie(S)$ is a rational surface, namely the involution associated with $\fie$ is
rational. The statement now follows from Theorem \ref{struttura}.
\end{proof}
\section{The examples}\label{sezexamples}
In this section we describe  all the  triples $(F, C,  G)$ that   occur
in Theorem \ref{struttura}.

We keep the notation introduced in \S 2 (cf. in particular Theorem
\ref{struttura}) and we introduce some more. We
denote by
$p_1\colon S\to F/G=\pp^1$ and $p_2\colon S\to C/G=\pp^1$  the two
isotrivial
pencils of
$S$. The singular fibres of $p_1$ and $p_2$  are multiples of smooth
curves, since the action of $G$ on $C\times F$ is free. We have seen in
the proof of Theorem \ref{struttura} that the multiple fibres of $p_2$
are 6 double fibres if $\si$ has 12 isolated fixed points and  5 double fibres
if
$\si$ has 10 isolated fixed points. We analyze now the  multiple fibres of $p_1$ and the
fixed set of
$\si$. Consider
$P\in F$ and denote by
$[P]\in F/G=\pp^1$ the equivalence class of $P$. The fibre of $p_1$ over
$[P]$ has multiplicity equal to the order of the stabilizer of $P$ in $G$.
The inverse image on $F\times C$ of the fixed locus of $\si$ consists of the
pairs $(P,Q)$ such that there exists $g\in G$ satisfying  $(\tau P, Q)=(gP,
gQ)$. When $g$ is the identity we get the divisorial part $R$ of the fixed
locus of $\si$.  The isolated fixed points of
$\si$ correspond  to solutions of the above equation with $P$  not a
fixed point
of $\tau$. Notice that in this case $g\inv \tau$ is a nontrivial
element of $\Aut(F)$ that fixes
$P$, hence if $\tau\in G$ then  the isolated fixed points of $\si$  lie on the multiple
fibres of
$p_1$. We denote by
$q_1, q_2\colon
\Si\to\pp^1$ the pencils induced by $p_1$ and
$p_2$. The general fibre of $q_2$ is a smooth rational curve, since $\si$ is
induced by the involution $\tau\times Id$ on $F\times C$.  The
inverse image in
$S$ of a general fibre of $q_1$ is connected iff the hyperelliptic
involution
$\tau$ is in
$G$.
\bigskip

In order to  describe the examples,  we need to study the possible group
actions on an hyperelliptic curve $F$. Given the group $G$, the action of $G$
on $F$ descends to an action on
$\pp^1=F/\tau$ that preserves the branch locus $\Delta$ of the
hyperelliptic
double cover $f\colon F\to \pp^1$. We denote by $\bar{G}$ the image of $G$
in
$\Aut(\pp^1)$.
 By identifying
 the unit sphere of $\R^3$  with $\pp^1$ via 
stereographic projection, one obtains an injection  $SO(3)\to
\Aut(\pp^1)$. In particular, the finite subgroups of $SO(3)$ give finite
subgroups of $\Aut(\pp^1)$ and, up to conjugacy, every finite subgroup of
$\Aut(\pp^1)$ arises in this way (cf. \cite{sfera}, Chapter III, or
\cite{enriques}, Libro Secondo, Capitolo I, \S 10). The finite subgroups of
$SO(3)$, up to conjugacy,
are the following: 
\begin{itemize}
\item[--] the cyclic groups $\Z_n$, generated by the rotation of $\frac{2\pi}{n}$ around the
$z-$axis;
\item[--] 
the  dihedral groups $D_n$, $n\ge 2$ (we set $D_2=\Z_2^2$), generated
by the rotation of $\frac{2\pi }{n}$ around the $z-$axis and by the reflection
in the line  $y=z=0$;
\item[--] the group 
   of symmetries of the tetrahedron, which is isomorphic to $A_4$;
\item[--] the group of symmetries of the cube (and of the octahedron),
which is isomorphic to 
$S_4$;
\item[--] the group of symmetries of the dodecahedron  (and of the
icosahedron), which is isomorphic to 
$A_5$.
\end{itemize}
Let $d=g(F)+1$ and let
$p(x_0,x_1)$ be a homogeneous polynomial of degree
$2d$ whose zero locus in $\pp^1$ is $\Delta$. The curve $F$ is isomorphic to the
curve of
equation $y^2=p(x_0,x_1)$ in the weighted projective plane $\pp(1,1,d)$. The
hyperelliptic involution
$\tau$ is
the restriction to $F$ of the automorphism of $\pp(1,1,d)$ given by
$(x_0,x_1,y)\mapsto (x_0,x_1,-y)$ and the hyperelliptic double cover $f\colon
F\to\pp^1$ is the restriction of the projection  $\pp(1,1,d)\to \pp^1$
defined by $(x_0,x_1,y)\mapsto (x_0,x_1)$. 
If we denote by $G_0$ the subgroup of $\Aut(F)$ generated by $G$ and $\tau$, then we have
 a central extension:
\begin{equation}\label{extension}
0\to <\tau>\to G_0\to\bar{G}\to 1.
\end{equation}
 Notice that  if $\tau\in G$ then  $G=
G_0$, while if $\tau\notin G$ then
$G$ is mapped isomorphically onto $\bar{G}$.

Consider now the quotient map  $\phi\colon SL(2,\C)\to PGL(1,\C)$ and set
$H:=\phi\inv \bar{G}$,  so that we have a surjective map
$H\to\bar{G}$ whose kernel is the subgroup generated by $-Id$. The group
$H$ acts on the space of homogeneous polynomials of degree $2d$ and
$p(x_0,x_1)$ is an eigenvector for this action. So there is a homomorphism
$\lambda\colon H\to \C^*$ defined by
$(h\inv)^*p(x_0,x_1)=\lambda(h)p(x_0,x_1)$. Since the degree of $p(x_0,x_1)$ is
even,  $-Id\in H$ is in the kernel of $\lambda$ and we have actually defined a
character
$\lambda\colon \bar{G}\to \C^*$. From now on we assume that $\lambda$ is  trivial 
 (this will always be the case in our  examples). In this case one can define a
homomorphism $H\to G_0$ by mapping $h\in H$ the automorphism of $F$   defined by
$(x_0,x_1,y)\mapsto (h(x_0,x_1),y)$. Since $d$ is even (by Theorem
\ref{struttura} either $d=4$ or $d=6$),
$-Id$ is in the kernel of
$H\to G_0$  and we have actually defined a homomorphism
$\psi\colon \bar{G}\to G_0$ that splits  the central extension
(\ref{extension}). We call $\psi$ the {\em canonical  splitting} of the sequence
(\ref{extension}).

Next we examine the fixed locus of the elements of $G_0$. We
may write
$G_0=\Z_2\times \bar{G}$,  identifying   $(1,0)$ with $\tau$ and
$\bar{G}$ with the image of $\psi$. Every  element $h$ of $G$ has two fixed
points
$P_1$, $P_2$ on $\pp^1$. If, say, $P_1\in \Delta$ then both $(0,h)$ and
$(1,h)$ fix the inverse image of $P_1$ in $F$. So assume that $P_1,P_2\notin
\Delta$ and write $f\inv(P_i)=\{Q_i^1, Q_i^2\}$, $i=1,2$. For each $i$,
either
$(0,h)$ fixes the points 
$Q_i^1$ and $Q_i^2$ or it  exchanges them. So, if the order of $h$ is odd, then
$(0,h)$
fixes the $4$ points $Q_i^j$  and $(1,h)$ acts freely
on
$F$. If the order of $h$ is even,
 let $h_0$ be an element of $SL(2,\C)$ that represents $h$ and let
$\alpha,\alpha\inv$ be the eigenvalues of $h_0$. Let   $(a, b)$ be
homogeneous coordinates for $P_1$, so that $Q_1^1$, $Q_1^2$ have coordinates
$(a,b,\pm c)$ with $c\ne 0$. Then $(0,h)$ maps $Q_1^1$,
$Q_1^2$ to the points $(\alpha a,
\alpha b, \pm c)=(a,b,\pm \alpha^{-d} c)$. So
$(0,h)$ fixes $Q_1^1$ and $Q_1^2$ iff $\alpha^d=1$, and the same is true for
$Q_2^1$
and $Q_2^2$. We conclude that if $\alpha^d=1$ then $(0,h)$ has $4$ fixed
points and
$(1,h)$ acts freely, while if $\alpha^d=-1$ then $(0,h)$ acts freely and
$(1,h)$ has
$4$ fixed points.

We will use  the following simple remark:
\begin{lem}\label{orbit}For $n\ge 2$ let  the dihedral  group $D_n$ of order
$2n$ act  on $\pp^1$ as follows:    
$$(x_0, x_1)\overset{r}{\mapsto}(x_0, e^{\frac{2\pi i}{n}} x_1);\qquad (x_0,
x_1)\overset{s}{\mapsto}(x_1,
 x_0),$$
where 
$r\in D_n$ is a rotation of order
$n$ and $s\in D_n$ is  a reflection.

Let $p(x_0,x_1)\in\C[x_0,x_1]$ be a  nonzero homogeneous polynomial of degree $2n$. Then
the zero locus of $p(x_0,x_1)$  is   an orbit of $D_n$  of order
$2n$ iff:

   $p(x_0,x_1)
=\alpha_0(x_0^{2n}+x_1^{2n})+\alpha_1x_0^nx_1^n$, with $\al_0\ne 0$, 
$2\al_0\pm\al_1\ne 0$.
\end{lem}
\begin{proof}
It is easy to check that if $p(x_0,x_1)$ is homogeneous of degree $2n$ and it vanishes on an
orbit of order $2n$, then  it  is of the
required form. Conversely, the zero locus of a polynomial 
$p(x_0,x_1)=\alpha_0(x_0^{2n}+x_1^{2n})+\alpha_1x_0^nx_1^n$ is invariant for $D_n$, hence it
is a union of orbits. 
The points of
$\pp^1$ that have nontrivial stabilizer are: $P_0=(1,0)$,
$P_1=(0, 1)$ and the points $Q_k=(1, e^{\frac{j\pi i}{n}})$,
$j=0\dots 2n-1$. A polynomial
$p(x_0,x_1)=\alpha_0(x_0^{2n}+x_1^{2n})+\alpha_1x_0^nx_1^n$\  vanishes at $P_0,
P_1$   iff $\al_0=0$ and it vanishes at $Q_j$ iff $2\al_0+(-1)^j\al_1=0$.  A
straightforward computation shows that if $p(x_0,x_1)$ vanishes at one of these
points then it has a multiple root there. This remark completes the proof.
\end{proof}
\bigskip

We are now ready to list all the possible types of minimal double planes of
general type $(S,\si)$ with $p_g=0$ and $K^2=8$. By Theorem \ref{struttura} $S$ is a  free
quotient 
$(F\times C)/G$ of the type described in Proposition \ref{costruzione} 
\smallskip

\noindent{\bf Type Ia:} Here $G=\Z_2^3$, $\bar{G}=D_2$ $(=\Z_2^2)$, $g(F)=3$,
$g(C)=5$ (cf. 
\cite{mp},
Example 4.2). We let  $\bar{G}$  be the subgroup
generated by:
$$(x_0,x_1)\stackrel{e_1}{\mapsto}(x_0,-x_1);\qquad 
(x_0,x_1)\stackrel{e_2}{\mapsto}(x_1,x_0)$$ 
and we take
$$p(x_0,x_1)=x_0^8+\alpha x_0^6x_1^2+\beta x_0^4x_1^4+\alpha x_0^2x_1^6+x_1^8,\quad
\pm 2\alpha+\beta+2\ne 0.$$
 Arguing  as in the proof of Lemma \ref{orbit}, one shows that
the condition on $\al$, $\be$ is equivalent to the fact that the zero locus $\Delta$ of
$p(x_0,x_1)$ is the union of two orbits of $\bar{G}$ of order 4.  We let
$G=G_0$ be the inverse image of $\bar{G}$ in $\Aut(F)$.  The
character $\lambda\colon \bar{G}\to \C^*$ is trivial and thus, as explained
before, we have a canonical  isomorphism of $G$ with $\Z_2\times \bar{G}$. The element
$(1,0)\in G$ corresponds to the hyperelliptic involution, hence it has 8
fixed points. The elements $(0,e_1)$,  $(0,e_2)$, $(0,e_1+e_2)$  have 4
fixed points each.   The remaining nonzero  elements, that we denote by
$f_1,f_2,f_3$,  act freely and are a set of generators  of $G$. By Proposition
2.1 of
\cite{ritaabel}, if we fix  distinct points $P_1\dots P_6\in \pp^1$ there
exists a smooth connected $G-$cover $C\to\pp^1$ such that the image of the fixed set of
$f_i$ is $P_{2i-1}+P_{2i}$, $i=1,2,3$, and the remaining nonzero elements of $G$ act freely
on
$C$. The curve $C$ is smooth of genus 5, hence we have the required example.

 The pencil
$p_1$ has
$5$ double fibres. The involution
$\si$ fixes pointwise the two 
double fibres   over the images in $\pp^1=F/G$  of the fixed
points of
$\tau$, and it  has
$4$ isolated fixed points on each of the remaining three.
\smallskip

\noindent{\bf Type  Ib:} Here $\bar{G}=D_4$, $G=\Z_2\times D_4$, $g(F)=3$,
$g(C)=9$. We let  $D_4$ act on   $\pp^1$ as in Lemma \ref{orbit} and we take: 
 $$p(x_0,x_1)=x_0^8+\alpha
x_0^4x_1^4+x_1^8,\quad\alpha\pm 2\ne 0.$$
By Lemma \ref{orbit}  the zero locus $\Delta$ of $p(x_0,x_1)$  is an orbit of  $D_4$ of
order 8.
 The 
 character
$\lambda\colon D_4\to \C^*$ is  trivial and  we consider the   decomposition
$G_0=\Z_2\times D_4$ given by the canonical splitting of the extension
(\ref{extension}).

 We now
examine the fixed points of the elements of $G=G_0$. The fixed set of the
hyperelliptic involution
$\tau=(1,0)$ consists of the $8$ Weierstrass points of $F$. By the previous 
discussion, the    elements $(0,sr^i)$, $i=0\ldots 3$,  have $4$
 fixed points each and the elements $(1,sr^i)$  act freely on
$F$. In addition, $(0, r^j)$ acts freely for $j$ odd and has $4$ fixed
points
for $j$ even, while $(1, r^j)$ has $4$ fixed points
for $j$ odd and acts freely for $j$ even, $j=0\ldots 3$.

Next we construct the curve $C$. Let $E$ be an elliptic curve and let
$\eta\in E$ be a point of order $4$. We define an action of $D_4$ on $E$ by:
$$z\stackrel{r}{\mapsto}z+\eta;\quad z\stackrel{s}{\mapsto}-z.$$ Denote by
$q\colon E\to E/D_4=\pp^1$ the quotient map.  Let
$h\colon
\pp^1\to\pp^1$ be a degree $2$ cover branched over two points that are not
branch points of $q$.   Let
$C$ be the fibre product of
$q$ and $h$, so that there is a commutative diagram:
\begin{equation}
\begin{CD}
C @>\tilde{h}>>E\\
@V\tilde{q}VV @V q VV\\
\pp^1@>h>>\pp^1
\end{CD}
\end{equation}

  To compute the genus of $C$ we apply the Hurwitz formula
to the double cover
$\tilde{h}\colon C\to E$, which is branched over $16$ points,  and we get
$g(C)=9$.  The curve  $C$ has a natural $\Z_2\times
D_4-$action as  a  fibred product.
The elements that have fixed points for this action are the following:
$(1,0)$ ($16$ fixed points),  $(0,sr^i)$,  $i=0\ldots 3$ ($8$ fixed points
each). In order to get a
free
action of $\Z_2\times D_4$ on
$F\times C$, we modify the action on $C$  by composing it with an
automorphism
$\psi$ of
$\Z_2\times D_4$ that maps the elements $(0,sr^i)$ to elements of the
form $(1, sr^j)$ and  exchanges
$(1,0)$ and $(1,r^2)$. For instance, one can take $\psi$ to be the automorphism
defined by $\psi(1,0)=(1,r^2)$, $\psi(0,r)=(0,r)$ and $\psi(0,s)=(1,s)$.
 With this action, the elements of $G$ that do not act
freely on
$C$ are:
$(1,r^2)$ ($16$ fixed points) and the elements $(1, sr^i)$, $i=0\ldots 3$
($8$
fixed points each). The multiple fibres of
$p_1$ are:  the double fibre over the image in $\pp^1=F/G$ of the fixed
points of $\tau$,  on which
$\sigma$ acts as the identity, a fibre of multiplicity
$4$ that contains $4$ isolated fixed points of $\si$, $2$ double fibres that contain
$4$
isolated fixed points of $\si$ each.
\smallskip

\noindent{\bf Type Ic:} Here $G\cong \bar{G}=S_4$, $g(F)=3$, $g(C)=13$. Up to
conjugacy, the action of $S_4$ on $\pp^1$ corresponds via the inverse of 
stereographic projection to the subgroup of rotations of the unit sphere in $\R^3$ that
preserve the cube of vertices $\frac{1}{\sqrt{3}}(\pm 1, \pm 1,\pm 1)$.  We take:
$$p(x_0,x_1)=x_0^8+(\alpha^4+\alpha^{-4})x_0^4x_1^4+x_1^8\subset \pp(1,1,4), \quad
\alpha=\frac{\sqrt{2}}{\sqrt{3}-1}.$$
 Namely, $F$ is the double cover
of
$\pp^1$  branched over the points corresponding to the vertices of the cube.
Since $S_4$ is not properly contained in
any finite subgroup of $PGL(1,\C)$,  in this case the central extension
(\ref{extension}) can be rewritten as:
$$0\to<\tau>\to
\Aut(F)\to S_4\to 1.$$

We claim that the character $\lambda\colon S_4\to \C^*$ is trivial.
Indeed, it is enough to check that there is
a transposition $\delta\in S_4$ with $\lambda(\delta)=1$. We take
$\delta$ to be defined by $(x_0,x_1)\mapsto (i x_1,i x_0)$.  The canonical
splitting of (\ref{extension}) gives a decomposition $\Aut(F)=\Z_2\times
S_4$. Next  we examine the fixed points of $\Aut(F)$. The
fixed points of $\tau=(1,0)$ are the $8$ branch points of the hyperelliptic
cover
$f\colon F\to\pp^1$. The element of $\Aut(F)$ that acts on $\pp(1,1,4)$ by
$(x_0,x_1,y)\mapsto (e^{\frac{ 3\pi i}{4}}x_0,e^{\frac{ 5\pi i}{4}}x_1, y)$ is the image
via $\psi$ of a cyclic permutation of order $4$ of $S_4$ and  it acts freely
on $F$. Since  the elements of order
$4$  of $S_4$ are all conjugated, this shows that  for every
$\delta$ of order $4$ the element  $(0,\delta)$ acts freely on $F$, 
while
$(1,\delta)$ has $4$ fixed points.
If $\delta$ is the transposition defined by $(x_0,x_1)\mapsto (i x_1,i x_0)$, then
$(0,
\delta)$ has $4$ fixed points and $(1,\delta)$ acts freely on $C$. Since all the
transpositions are conjugated in $S_4$, we conclude that the same is true
for
every transposition
$\delta$. If $\delta$ is a $3-$cycle then both $(0,\delta)$ and $(1,\delta)$
have
$2$ fixed points.  We let
$G\cong S_4$  be the subgroup
$\{(\epsi(\delta),
\delta)|\delta \in S_4\}\subset \Z_2\times S_4$, where $\epsi(\delta)$ is the sign
of the permutation $\delta$.
 Now  consider a generic degree
$4$ map
$h\colon
\pp^1\to\pp^1$. The branch locus of $h$ consists of $6$ distinct points,
each
corresponding to a simple ramification point. The Galois closure
$C\to\pp^1$ of $h\colon \pp^1\to\pp^1$ is an
$S_4-$cover branched on $6$ points such that the  only elements of $S_4$ that
 do not act freely on $C$ are the transpositions. The Hurwitz formula gives
$g(C)=13$ and the diagonal action of $G$ on $F\times C$ is free.
Notice that in this case the hyperelliptic involution $\tau$ is not in $G$,
and
thus the pullback on $S$ of a general fibre of $q_1$ is disconnected.
The multiple fibres of $p_1$ are: a triple fibre, that is fixed pointwise by
$\si$, and two fibres of multiplicity $4$ that are interchanged by $\si$.
The
$12$ isolated fixed points of $\si$ all lie on the  smooth fibre of $p_1$
over
the image in $F/G$ of the fixed points of the elements $(0,\delta)\in G$,
with
$\delta$ a transposition.
\smallskip

\noindent{\bf Type  Id:}  Here $G=\Z_2\times S_4$, $g(F)=3$, $g(C)=25$.
The
curve
$F$ is the same as in type Ic.  We take $G=\Aut(F)$ and we
consider again 
the   direct sum decomposition $G=
\Z_2\times S_4$ given by the canonical splitting of (\ref{extension}).
Let $p\colon E\to \pp^1$ be a degree $4$ cover with the following
properties:
\begin{itemize}
\item[a)] $E$ is
a smooth curve of genus $1$;

\item[b)] the branch locus of $p$ consists of points
$P_1\ldots P_4, Q_1,Q_2$ such that $p$ has a simple ramification point over
$P_1\ldots P_4$ and  it has two simple ramification points over $Q_1,Q_2$;

\item[c)] the Galois closure  $q\colon D\to\pp^1$ of $p$ has Galois group equal
to
$S_4$.
\end{itemize}
We remark that the existence of such a cover $p$ can be shown by using the
classical Riemann construction. Let $h\colon \pp^1\to\pp^1$ be the double
cover branched on $Q_1, Q_2$ and let $C$ be the normalization of the fibre
product  of
$h$ and
$q$. We have a commutative diagram:
\begin{equation}
\begin{CD}
C @>\tilde{h}>>D\\
@V\tilde{q}VV @V q VV\\
\pp^1@>h>>\pp^1
\end{CD}
\end{equation}
Applying the Hurwitz formula to $q$ one obtains $g(D)=13$. The map
$\tilde{h}$
is
\'etale of degree
$2$.  The group actions on
$\pp^1$ and
$D$ induce an action
of
$\Z_2\times S_4$ on $C$ such that $C/(\Z_2\times\{0\})=D$ and $C/(\{0\}\times
S_4)=\pp^1$.   The elements with nonempty fixed locus are those of the form
$(0,\delta)$,  with $\delta$ a transposition,  and
$(1,\delta)$, with
$\delta$ a double cycle. Each of these elements fixes $16$ points.
We claim that the curve $C$  is connected. Assume otherwise. Then 
$C$ is the disjoint union of two connected components, isomorphic to $D$, that are
exchanged by $(1,0)$. The subgroup $G$ of the elements that do not exchange the
components  of $C$ is $\{0\}\times S_4$ and it should contain all the elements
that do not act freely on $C$, contradicting the above analysis of the fixed
points of the elements of
$G$. The Hurwitz formula then gives 
$g(C)=25$.

In order to get a free action on $F\times C$,  we have  modify the
given
action on $C$ by composing it with an automorphism $\psi$ of $G$ of the form  $(a,
\delta)
\mapsto (a+\epsi(\delta), \psi_1(\delta))$, where $\psi_1$ is an automorphism of
$S_4$ and $\epsi(\delta)$ is the sign of $\delta$. For instance one can take $\psi_1=Id$.

The multiple fibres of $p_1$ are: a $6-$tuple fibre, which is fixed by $\si$
pointwise, a double fibre containing $8$ isolated fixed points of $\si$
and a
$4-$tuple fibre containing $4$ isolated fixed points of $\si$.
\smallskip

\noindent{\bf Type II:} This  is similar to type  Ic. Here
$g(F)=5$,
$g(C)=16$ and
$G\cong\bar{G}=A_5$. The group $A_5$ acts on the unit sphere in $\R^3$ as the group of
rotations that preserve the dodecahedron (and the icosahedron).  Identifying
as
usual the sphere with $\pp^1$
via stereographic projection, we get an action of
$A_5$
on
$\pp^1$. The special orbits are: one orbit of order $12$, one orbit of order
$30$ and one orbit of order $20$. Let $p(x_0,x_1)$ be a homogeneous polynomial
of
degree
$12$ whose zeros on $\pp^1$ are the points of the orbit of order $12$ and
let
$F\subset \pp(1,1,6)$ be the hyperelliptic curve of genus $5$ defined by
$y^2=p(x_0,x_1)$.   Since $A_5$ is not
properly
contained in any finite subgroup of $\Aut(\pp^1)$, in this case the  central extension
(\ref{extension}) can be rewritten:
$$0\to<\tau>\to \Aut(F)\to A_5\to 1.$$ As we have already remarked, the
character
$\lambda\colon A_5\to \C^*$ is necessarily trivial.  Thus the  above sequence
is split  and one has  a  decomposition
$\Aut(F)=\Z_2\times A_5$, which is unique. Now we examine the fixed points
 of the elements of 
$\Aut(F)$.  The hyperelliptic involution $\tau=(1,0)$ has $12$ fixed points. If
$\delta$
is a
$5-$cycle in
$A_5$, then both
$(0,\delta)$ and
$(1,\delta)$ have $2$ fixed points. If $\delta$ is a $3-$cycle, then
$(0,\delta)$
has $4$ fixed points and $(1,\delta)$ acts freely. Now, applying the Hurwitz
formula to the covering $F\to
\pp^1=F/A_5$,  one sees that the remaining elements of $A_5$ act freely. It
follows that   the elements $(1,\delta)$, with $\delta$
a
double cycle, have $4$ fixed points each.
Now let $C\to\pp^1$ be a Galois cover with Galois group $A_5$, branched over
$5$ points and such that the only elements that do not act freely on $C$
are the double cycles (the existence of such a cover can be shown using the
Riemann construction). Then the genus of $C$ is $16$ by the Hurwitz formula
and the diagonal action of $G=A_5$ on $F\times C$ is free, as required.
As in type  Ic, the hyperelliptic involution $\tau$ is not in
$G$, and thus the pullback on $S$ of a general fibre of $q_1\colon \Si\to\pp^1$ is
disconnected.
The  multiple fibres of $p_1\colon S\to\pp^1$ are: a $5-$tuple fibre that is
fixed by $\si$ pointwise and $2$ triple fibres that are exchanged by $\si$.
The $10$ isolated fixed points all belong to the  fibre of $p_1$ over the
image
of a point $P\in F$ whose stabilizer in $\Aut(F)$ is generated by an element
$(1,\delta)$ with $\delta$ of order $2$.
\bigskip
\section{The classification}\label{sezclassifica}
We are now ready to state and prove  our main result. We keep all the notation of the
previous sections.
\begin{thm}\label{classifica}Let $(S,\si)$ be a minimal double plane of general
type with
$p_g(S)=0$ and
$K^2_S=8$. Then
 $(S,\si)$ belongs to one of the types Ia, Ib, Ic, Id, II
described in \S {\rm \ref{sezexamples}}.
\end{thm}

 We say that  a minimal double plane of general type $(S,\si)$ with
$K^2=8$ and  $p_g=0$ is of  type I if it  belongs to one of the types Ia,
Ib, Ic, Id.   By Theorem \ref{struttura} and Theorem \ref{classifica}, the double
planes of type I are characterized by the fact that they have 12 isolated fixed
points and the double planes of type II
 by the fact that they have 10 isolated fixed points.

Before proving Theorem \ref{classifica} we remark the following:
\begin{cor}\label{unique}
Let $(S,\si)$ be a minimal double plane of general type with $K^2_S=8$ and
$p_g(S)=0$. Then:
\begin{itemize}
\item[i)] $\si$ is the only rational involution of $S$;
\item[ii)]$(S,\si)$ belongs to exactly  one of the types Ia, \dots Id, II;
\item[iii)]the bicanonical map
of
$S$ is birational iff $(S,\si)$ is of type II.
\end{itemize}
\end{cor}
\begin{proof}
Assume by contradiction that $S$ has two distinct  rational involutions $\si_1$
and
$\si_2$. By Theorem \ref{struttura}, to each involution $\si_i$ there
corresponds an isomorphism of 
$S$ with a quotient
$(F_i\times C_i)/G_i$, where $F_i$ is hyperelliptic and $\si_i$ is induced by
$\tau_i\times Id$,  $\tau_i$ being the hyperelliptic involution of $F_i$.
The projections of $F_i\times C_i$  onto the two factors induce pencils $p_1^i$,
$p_2^i\colon S\to\pp^1$, $i=1,2$. A base point free pencil of $S$ is determined
uniquely by the span of the cohomology class of one of its fibres, which is an isotropic
subspace  of $H^2(S, \Q)$. Since $h^2(S, \Q)=2$, $S$ has exactly two free pencils and thus 
 we have $\{p^1_1,p^1_2\}=\{p^2_1,p^2_2\}$.
Since the involutions $\si_1$ and $\si_2$ are different, we have $p^1_1=p^2_2$,
$p^2_1=p^1_2$, or, equivalently, $C_1=F_2$ and  $C_2=F_1$. Looking at the pair  $(g(F),
g(C))$  for the types listed in   Theorem
\ref{classifica} one sees that this cannot happen. This proves i). Statement ii) follows
from i) and from Theorem \ref{classifica} by the same argument, since  the   pair
$(g(F),g(C))$ is different for each type of double planes.   Now statement iii) follows
from  i), ii) and  Corollary
\ref{bica}.
\end{proof}

\begin{proof}[Proof of Theorem \ref{classifica}] 
The proof is long and it is divided in several steps.
\vspace{0.1cm}

By Theorem \ref{struttura}
there exist a curve $C$, an hyperelliptic curve $F$ of genus 3 or 5  and a
group $G$ that acts faithfully on $C$ and $F$ such that
$S=(F\times C)/G$ and
$\si$ is induced by the involution $\tau\times Id$ of $F\times C$, where
$\tau$ is the hyperelliptic involution of $F$.    
We denote by $A\subset G$ (resp. $B\subset G$) the subset of elements
$\ne 1$ that do not act freely on $F$ (resp. $C$).
\begin{step}\label{condizioni}
The following conditions  are satisfied:
\begin{itemize}

\item[a)] $A\cap B$ is empty;

\item[b)] $A$ and $B$ are a union of conjugacy classes of $G$;
\item[c)] both $A$ and $B$ generate $G$ (in particular both $A$ and $B$ are
nonempty);
\item[d)] the elements of $B$ have order $2$.
\end{itemize}
\end{step}
Condition  a) is equivalent to the fact that the diagonal
action of $G$ on 
$F\times C$ is free (cf. Theorem \ref{struttura}, a)). Condition b) follows from the
definition of $A$ and $B$. Condition c) follows from the fact that by Theorem
\ref{struttura}, b) the curves  $F/G$ and
$C/G$ are rational. In fact, if $H$ is the subgroup of $G$
generated by $B$, then $C/H\to
C/G=\pp^1$ is a connected \'etale cover, hence $H=G$.  The same argument shows that $G$
is generated by $A$.
To prove    d), we recall that for a point  $P\in C$ the multiplicity  of the fibre of
$p_2\colon S\to C/G$ over 
$[P]\in C/G$ is equal to the order of ramification at $P$ of the quotient map
$C\to C/G$, which in turn is equal to the order of the stabilizer  of
$P$ in $G$. Since all the  multiple fibres of $p$ are double (cf. proof of Theorem
\ref{struttura}) $B$ consists   of elements of order 2.

\begin{step}\label{diedrale}
$G$ is not one of the following groups:
$$\Z_n,\ n\ge 2, \qquad
D_n, \ n\ge 2$$
\end{step}
\noindent By Step \ref{condizioni}, c), d), $G$ is generated by elements of
order
$2$. So if $G$ is cyclic,
then it is equal to $\Z_2$, but then  $G$ does not have two disjoint sets of
generators, contradicting Step \ref{condizioni}, a), c). The same argument rules out
$D_2=\Z_2^2$. So assume that
$G$
is a dihedral group $D_n$, $n\ge 3$.  Looking at the conjugacy classes of the
elements of order 2 of $D_n$, one sees that  conditions  b),  c) and d)  of Step
\ref{condizioni} imply that $B$   contains all the reflections. Then $A\subset
D_n\setminus B$ does not generate
$D_n$, a contradiction to Step \ref{condizioni}, c).
\bigskip

As in \S \ref{sezexamples}, we denote by $\bar{G}$ the finite subgroup of automorphisms of
$\pp^1=F/\tau$ induced by $G$.
\begin{step}\label{Gbar}The group  $\bar{G}$  is one of the   following:
$$ D_n, \ n\ge 2;\quad S_4;\quad A_5.$$
\end{step}
\noindent This follows by Step
\ref{diedrale} and by the classification of the finite subgroups of
$PGL(1,\C)$ that we have recalled at the beginning of  \S \ref{sezexamples} ($A_4$ is
excluded since it is not generated by elements of order $2$, contradicting Step
\ref{condizioni}, d)).
\bigskip

 We denote by
$\bar{A}$, $\bar{B}$ the images of
$A$, $B$ in $\bar{G}$. We notice that by Step \ref{condizioni} $\bar{A}$ and $\bar{B}$ are
sets of generators of $\bar{G}$ and are stable under conjugacy. The 
elements
of $\bar{B}$ have order $2$.
\begin{step}\label{free}
If $\bar{G}=D_n$, $n\ge 2$, then $\bar{G}$ acts freely on the branch locus
$\Delta$ of
$F\to F/\tau$.
\end{step}
\noindent By Step \ref{diedrale}, if $\bar{G}=D_n$ then $G$ is not isomorphic to $\bar{G}$
and thus the exact sequence
(\ref{extension})  becomes:
$$0\to<\tau>\to G\to D_n\to 1.$$
 Arguing as in the proof of Step
\ref{diedrale}, one sees that
$\bar{B}$ contains all the reflections of $D_n$. By Step \ref{condizioni}, a), d), this is
the same as saying that every reflection can be lifted to an element of order 2 of $G$ 
that acts freely on $F$. 
On the other hand, if $\bar{P}\in \Delta$ is fixed by an element $\bar{g}\in \bar{G}$, 
then  the inverse image
$P\in F$ of
$\bar{P}$ is fixed by both the elements $g_1,g_2\in G$ that lift
$\bar{g}$.
Hence the reflections of
$\bar{G}=D_n$ act freely on $\Delta$. Assume now that
$\bar{P}\in \Delta$ is fixed by a rotation $r$ of $D_n$. We may assume that $r$ has order
$n$. Let
$r_1, r_2=\tau r_1\in G$ be the  lifts of
$r$. The stabilizer $H$ of $P$ in $G$ is cyclic, since $G$ acts faithfully 
on the smooth
curve
$F$,  and it is generated by 
$r_1$,
$r_2$ and
$\tau$. Hence, up to exchanging $r_1$ and $r_2$,   we may
assume that $r_1$ generates $H$ and, in particular, 
$\tau=r_1^n$.  Let $s\in D_n$ be a reflection and let $s_1\in G$ be an element of
order 2 that lifts $s$. Since $srs=r\inv$,  we 
either have
$s_1r_1s_1=r_1\inv$ or $s_1r_1s_1=\tau r_1\inv=r_1^{n-1}$. In the former case, 
$G$ is isomorphic to
$D_{2n}$, contradicting Step \ref{diedrale}. The latter case can occur only for
$n$ even, since for $n$ odd  $r^{n-1}$ does not generate $H$.  In
this case the elements
$sr_1^i$ and $sr_1^i\tau$  have order 4 for $i$ odd and thus they  are not contained in $B$,
contradicting the fact that $\bar{B}$ contains all the reflections.
\bigskip

\begin{step}\label{even}
If $\bar{G}=D_n$, then $n$ is even.
\end{step}
\noindent
By Step \ref{diedrale}, if $\bar{G}=D_n$ then $G$ and  $\bar{G}$ are not isomorphic and  the
exact sequence
(\ref{extension}) gives:
$$0\to<\tau>\to G\to D_n\to 1.$$
 Assume by contradiction
that
$n$ is odd. Then a rotation $r\in D_n$ of order $n$ can be lifted to $r'\in G$ of
order
$2n$. Arguing as in the proof of  Step \ref{free} one shows that 
$G$ is dihedral, contradicting Step \ref{diedrale}.
\bigskip

As we have already remarked, the curve $F$ is defined by an equation of the
form $y^2=p(x_0, x_1)$ in the weighted projective plane  $\pp(1,1,d)$, where
$d=g(F)+1$ and $p$ is a homogeneous polynomial of degree $2d$. 
\begin{step}\label{cases}Up to isomorphism, there are the following  possibilities for
$\bar{G}$,
$g(F)$ and  
$p(x_0,x_1)$:
\begin{itemize}
\item[i)] $\bar{G}=D_2$, $g(F)=3$, $p(x_0,x_1)=x_0^8+\alpha x_0^6x_1^2+\beta
x_0^4x_1^4+\alpha x_0^2x_1^6+x_1^8$,
$2\pm 2\alpha+\beta\ne 0$;
\item[ii)] $\bar{G}=D_2$, $g(F)=5$, $p(x_0,x_1)=x_0^{12}+\alpha
x_0^{10}x_1^2+\beta x_0^8x_1^4+\ga x_0^6x_1^6+ \beta x_0^4x_1^8 +\alpha
x_0^2x_1^{10}+x_1^{12}$,
$2+ 2\alpha+2\beta+\ga\ne 0$, $2- 2\alpha-2\beta+\ga\ne 0$;
\item[iii)]  $\bar{G}=D_4$, $g(F)=3$, $p(x_0,x_1)=x_0^8+\alpha
x_0^4x_1^4+x_1^8$,
$2\pm\alpha\ne 0$;
\item[iv)] $\bar{G}=D_6$, $g(F)=5$, $p(x_0,x_1)=x_0^{12}+\alpha
x_0^6x_1^6+x_1^{12}$,
$2\pm\alpha\ne 0$;
\item[v)] $\bar{G}=S_4$, $g(F)=3$, 
$p(x_0,x_1)=x_0^8+(\alpha^4+\alpha^{-4})x_0^4x_1^4+x_1^8\subset \pp(1,1,4)$
where $\alpha=\frac{\sqrt{2}}{\sqrt{3}-1}$ (the zeroes of $p$ are the points
of the orbit of order 8 of $S_4$);
\item[vi)] $\bar{G}=A_5$, $g(F)=5$, $p(x_0,x_1)$ is a polynomial of degree 12
whose zeroes are the elements of the orbit   of order 12 of $A_5$.
\end{itemize}
In all cases,  the exact sequence (\ref{extension}) is split.
\end{step}
We recall that the zero set $\Delta$ of $p(x_0,x_1)$ is the branch
locus of  the hyperelliptic double cover $f\colon F\to \pp^1$ and it is invariant under
the action of
$\bar{G}$. If
$\bar{G}=A_5$, then  we have necessarily case vi), since the smallest orbit
of $A_5$ on
$\pp^1$ has 12 elements and there is only one such orbit.
If $\bar{G}=S_4$, then a similar argument shows that we have case v). 
By Step \ref{Gbar}, the only remaining possibility is that $\bar{G}$ is a
dihedral group $D_n$. If this is the case, then Step \ref{free} and   Step
\ref{even}  imply that the only possibilities for $\bar{G}$, and $g(F)$ are i)--iv).
The fact that $p(x_0,x_1)$ is of the form stated above follows from Lemma \ref{orbit} in
cases iii), iv), it has been proven while describing type Ib for case i) and it can be
proven in the same way for case ii). The fact that the exact sequence  (\ref{extension})
is split has already been checked in \S \ref{sezexamples} while describing the
various types of double planes except for case ii), where one can use exactly the
same argument.
\bigskip

\begin{step} If  $\bar{G}=D_2$, then $(S,\si)$ is of type Ia.
\end{step} 
Up to a suitable choice of homogeneous coordinates, we may assume
that the action of $\bar{G}$  on $\pp^1$ is  the one described  for surfaces of
type Ia. By Step
\ref{diedrale}, the groups $G$ and $\bar{G}$ are not isomorphic, hence
(\ref{extension}) gives the exact sequence:
$$0\to<\tau>\to G\to  D_2\to 1$$
which is split by Step \ref{cases}.  The curve $F$ is as in case i) or ii) of Step
\ref{cases} and, with respect to the canonical splitting  $G\cong
\Z_2\times D_2$, each of the nonzero elements of
$\{0\}\times D_2$ has 4 fixed points on $F$. Thus $A$ contains $\tau$ and three  
elements
$e_1,e_2,e_3$ that map to the three nontrivial elements of $D_2$. By conditions a) and c) of
Step
\ref{condizioni}, 
$B$ contains the remaining nonzero elements $f_1,f_2,f_3$ of $G$,  which are a set of
generators.
  We consider the quotient
map $h\colon C\to C/G=\pp^1$. By Theorem
$2.1$ of
\cite{ritaabel}, the image via $h$ of the fixed set of $f_i$ has cardinality
divisible by $2$ for  $i=1,2,3$.  Then the Hurwitz formula implies
$g(C)\equiv 1$
$(\!\!\!\mod 4)$. Hence case ii) of Step \ref{cases} is excluded and we have case
i) of Step
\ref{cases}.  Now, using again Theorem 2.1 of \cite{ritaabel}, it is easy to check that we
have a double plane of type   Ia.
\bigskip

\begin{step} If $\bar{G}=D_n$, $n\ge 3$, then $(S,\si)$ is of type Ib.
\end{step}

Here we have either case iii) or case iv) of Step \ref{cases}. We will show that case iii)
corresponds to type Ib and that iv) does not occur.

We consider  case iii) first.  By Step
\ref{diedrale} the groups $G$ and $\bar{G}$ are not isomorphic, hence
(\ref{extension}) gives the exact sequence:
$$0\to<\tau>\to G\to  D_4\to 1$$
which is split by Step \ref{cases}. 
The curve $F$ is the same as in type Ib and $g(C)=9$.   Using the canonical
splitting
$G\cong
\Z_2\times D_4$  one sees that  $A$ contains $\tau=(1,0)$ and $(0,
sr^i)$, $i=0\dots 3$ (cf. the description of type Ib in \S \ref{sezexamples}). It
follows by Step
\ref{condizioni} that $B=\{(1, r^2), (1,sr^i), i=0\dots 3\}$.
 Let $\psi\colon C\to
C/G=\pp^1$
be the projection onto the quotient and let $\nu$ (resp. $\mu$) be
the number of branch points of $\psi$ that are images of the fixed points of
$(1,r^2)$ (resp.  of the elements $(1,sr^i)$). So
$\nu,\mu>0$ and the Hurwitz formula gives $\nu+\mu=6$.
Set
$E:=C/(1,r^2)$ and consider the quotient map $C\to E$: the Hurwitz formula
gives
$16=2(2g(E)-2)+8\nu$ and thus $\nu\le 2$. Let $H$ be the
subgroup
generated by $(1,s)$ and $(0, r)$ and set $D:=C/H$. The subgroup $H$ is
isomorphic
to $D_4$ and the Hurwitz formula gives: $16=8(2g(D)-2)+8\mu$.
Since
$\mu=6-\nu\ge 4$,   we have:
$g(D)=0$,  $\mu=4$, $\nu=2$ and $g(E)=1$.  The group
$G/<(1,r^2)>$ is isomorphic to  $D_4$ and it  
acts on $E$ in such a way that only the  reflections have fixed points. It
follows
that the action on $E$ must be the one described for type Ib. In
addition, we
have a commutative diagram:
\begin{equation}
\begin{CD}
C @>\tilde{h}>>E\\
@V\tilde{q}VV @V q VV\\
C/H=\pp^1@>h>>E/D_4=\pp^1
\end{CD}
\end{equation}
The branch points of $h$ are the images of the points of $C$ fixed by $(1,r^2)$ and
the branch points of $q$ are the images of the points fixed by the elements of type
$(1, sr^i)$.  In particular,  $q$ and $h$ have no common branch point and thus
their fibre product is smooth and connected. It follows that $C$ is isomorphic to the fibre
product of $q$ and $h$ and that  this is type
 Ib.

Now we consider  case iv) of Step \ref{cases}.
By Step
\ref{diedrale} the  groups $G$ and $\bar{G}$ are not isomorphic, hence
(\ref{extension}) gives the exact sequence:
$$0\to<\tau>\to G\to  D_6\to 1$$
which is split by Step \ref{cases}. 
Arguing as in the description of type Ib, one sees that with respect to the canonical 
decomposition
$G=\Z_2\times D_6$  the elements of order $2$ that act freely on
$F$ are the
following:  $(0,sr^i)$, $i=0\dots 5$ and $(0,r^3)$. Since these elements do not
generate
$G$, this case does not occur by Step \ref{condizioni}.
\bigskip

\begin{step} In case v) of Step \ref{cases}, $(S,\si)$   is either of  type Ic
or of type Id.
\end{step}
In this case, the curve $F$ and its automorphism group have been analyzed in detail in the
description of surfaces of type Ic in \S
\ref{sezexamples}. In particular, it has been shown that, with respect to
the canonical  splitting $G=\Z_2\times S_4$,  the elements of order
$2$ that act freely on $F$
are those of the form  $(1,\delta)$ with $\delta$ of order $2$.
 Then, by Step \ref{condizioni},  we have
two
possibilities:
\begin{itemize}
\item[a)] $G= \Aut(F)\cong \Z_2\times S_4$ and $B$ consists of all the elements
$(1,\delta)$
with
$\delta$ of order $2$.
\item[b)]  $G=\{(\epsi(\delta), \delta)\}\cong S_4$, where $\epsi$ denotes the sign
of $\delta$, and $B$ consists of all the
transpositions.
\end{itemize}

We consider  case b) first. One has $g(C)=13$ and by the Hurwitz formula  the branch
locus of the
quotient map
$C\to C/S_4=\pp^1$ consists of $6$ points. Thus each transposition of $S_4$
fixes
$12$ points. Let $H$ be a subgroup of $S_4$ isomorphic to $S_3$. By the
Hurwitz
formula, the curve $C/H$ is rational. In addition, the induced map
$C/H=\pp^1\to
C/S_4=\pp^1$ is a degree $4$ morphism with $6$ simple branch points. This shows
that case
b) corresponds to type Ic

Now assume we are in case a) and denote by $\nu$ the number of branch
points that
are the images of the fixed points of the elements $(1,\delta)$ with
$\delta$ a
double cycle and by $\mu$ the number of branch points that are the images
of the
fixed points of the elements $(1,\delta)$ with $\delta$ a transposition.  So
$(1,\delta)$ has $8\nu$ fixed points if $\delta $ is a double cycle and
it has
$4\mu$ fixed points if $\delta$ is a transposition. We have
$\nu,\mu>0$ and
the Hurwitz formula gives again $\nu+\mu=6$. Denote by $H$ the subgroup
of $G$
generated by the elements
$(1,\delta)$ with $\delta$ a double cycle. $H$ is isomorphic to $\Z_2^3$ and
the
Hurwitz formula applied to the quotient map $C\to C'=C/H$ gives
$48=8(2g(C')-2)+24\nu$, which implies:  $g(C')=1$, $\nu=2$, $\mu=4$.
Denote
by
$P_1\dots P_4\in \pp^1$ (respectively $Q_1,Q_2$) the branch points corresponding to
the elements
$(1,\delta)$
with
$\delta$ a transposition (respectively  a double cycle). Denote by
$K$ the subgroup of
$G$ generated by the elements
$(1,\delta)$, where
$\delta$ is a transposition. $K$ is isomorphic to $S_4$ and the quotient
curve $C/K$ is
rational by the Hurwitz formula. The quotient map $h\colon C/K=\pp^1\to C/G=\pp^1$
is a
degree $2$ cover branched over $Q_1$ and $Q_2$. We set
$D:=C/\Z_2\times\{0\}$. The
quotient map $\tilde{h}\colon C\to D$ is an \'etale double cover and thus
$D$ has
genus $13$.
Summing up, we have a commutative diagram:
\begin{equation}
\begin{CD}
C @>\tilde{h}>>D\\
@V\tilde{q}VV @V q VV\\
\pp^1@>h>>\pp^1
\end{CD}
\end{equation}
Thus the  curve $C$ is obtained from $q$ and $h$ by base change and
normalization.  To show that this is type Id, we
consider
the action of
$G/\Z_2\times\{0\}\cong S_4$  on
$D$. The  quotient map
$q\colon D\to C/G=\pp^1$ is branched over
$P_1\dots P_4,Q_1,Q_2$. If $K_1$ is a subgroup of $S_4$ isomorphic to $S_3$
and we
write $E:=D/K_1$, then $E$ has genus $1$ and the induced map $E\to
C/G=\pp^1$ is a
degree $4$ cover with one simple ramification point above $P_1\dots P_4$ and
$2$
simple ramification points over $Q_1$ and $Q_2$, whose Galois closure is
$q\colon
D\to \pp^1$, as required.
\bigskip

\begin{step}In case vi) of Step \ref{cases}, $(S,\si)$ is of type II
\end{step}

As  explained in the description of type II,
there is
a unique splitting $\Aut(F)=\Z_2\times A_5$ and the elements of order $2$
that act
freely on $F$ are the elements $(0,\delta)$, where $\delta$ is a double
cycle of
$A_5$. The subgroup generated by these is $A_5$, thus we have $G=A_5$ and
$g(C)=16$.
Arguing as in the previous steps, it is easy to check that this is type II.  
\end{proof}
\section{Plane models}\label{sezmodels}

 A {\em plane model} of a double plane  $(S,\si)$  is a finite degree 2 morphism
$X\to\pp^2$ such that $X$ is a normal surface and there exists a commutative diagram:
$$\begin{CD}
S@>>>X\\
@VVV @VVV\\
\Si=S/\si@>>>\pp^2
\end{CD}$$
such that the horizontal arrows denote  birational maps. It is clear that a plane model is
by no means unique, since one can compose the map  $X\to\pp^2$ with a birational
transformation of $\pp^2$, solve the indeterminacy of the resulting map and consider
the Stein factorization,   thus obtaining a different plane model. 

In this section we compute plane models for the minimal double planes of general type with
$K^2=8$ and $p_g=0$ using the classification of \S  \ref{sezclassifica}. It turns out that
the double planes of type I can all be realized  using  a construction suggested by Du
Val in \cite{duval}
 and the double planes of type II can be  obtained in a similar  way.

 A plane model 
$X\to
\pp^2$ is determined by its branch curve
$B$, which is a reduced curve of even degree. The invariants 
$\chi$ and $K^2$  of the minimal resolution of $X$  can be computed in terms of  the degree
of
$B$ and of the type of its singularities, while the values of $p_g$ and $q$ depend also on
the mutual positions of the singularities of $B$ (cf. \cite{bpv}, Ch. III, \S 7).  In our
case
  the ``expected number'' of parameters for $B$ is negative and
the singularities of $B$ satisfy a general position property (cf. Theorem \ref{planeI}, iv)
and Theorem \ref{planeII}, iv)).  Therefore  it seems   very difficult to construct
these plane models directly, let alone classify them.

We call a singular point of a plane curve {\em simple} if it is solved by  blowing up the
singular point in the plane once. A point of type $(m,m)$ is a point of multiplicity $m$
with an infinitely near simple point also  of multiplicity $m$.

The following theorem describes  the Du Val plane model for type I.
\begin{thm}\label{planeI} Let $(S,\si)$ be a minimal double plane of general type with
$K^2=8$ and
$p_g=0$. 

If $(S,\si)$ is  of type I, then there
 is a plane
model
$X\to\pp^2$ with branch curve $B$ such that:
\begin{itemize}
\item[i)] $B=C_{16}+L_1+\dots + L_6$, where $C_{16}$ has degree 16 and $L_1\dots L_6$ are
distinct lines;
 
\item[ii)] the singularities of
$C_{16}$ are:
\begin{itemize}
\item[--] a  singular  point $P$  of multiplicity 8 such that   $P\in L_1\dots L_6$;

\item[--]  6 points $R_i$ of type $(4,4)$ such that $R_i\in L_i$ and $L_i$ is tangent to
$C_{16}$ at $R_i$,
$i=1\dots 6$;
\end{itemize}

\item[iii)] the singularities of $B$ are
solved by blowing up $P$, $R_1\dots R_6$ and the points $S_i$ infinitely near to $R_i$ in the
direction of $L_i$, $i=1\dots 6$;

\item[iv)] there is no conic containing  $R_1\dots R_6$.
\end{itemize}

One or two of the $R_i$ can be infinitely near to $P$. The curve $C_{16}$ has two
irreducible  components of degree  8 if $(S,\si)$ is of type Ia and it is irreducible
otherwise.

Conversely, given a plane  curve $B$ as in i)--iv), the double cover
$X\to\pp^2$ branched on $B$ is a plane model of a minimal double plane of general type with
$K^2=8$ and $p_g=0$, of type I.
\end{thm}
\begin{proof}Denote by $\tilde{S}$ the blow-up of $S$ at the
isolated fixed points of
$\si$,  by
$\tilde{\si}$ the involution of $\tilde{S}$ induced by $\si$, and by $Y$ the quotient
surface $\tilde{S}/\tilde{\si}$. $Y$ is the minimal resolution of $\Si=S/\si$.

We denote by $\F_e$ the geometrically  ruled surface
$\Proj(\OO_{\pp^1}\oplus\OO(e))$,
$e\ge 0$, by
$C_{\infty}$ the class of a section of $\F_e$ with self-intersection $-e$ and by $\Phi$ the
class of a ruling.  Arguing as in the proof of Theorem 4.4 of
\cite{nodi}, one can show that  $Y$ is obtained from  a surface $\F_e$ by blowing up $6$
pairs of points $P_i, P'_i$ such that:
\begin{itemize}
\item[a)]
$P_i, P_j$ lie on different rulings of $\F_e$ if $i\ne j$;
\item[b)] for every
$i=1\dots 6$ the point $P'_i$ is  infinitely near   to $P_i$  in the
direction of $\Phi$;
\item[c)] the pull back of the fibration $Y\to\F_e\to\pp^1$ is the hyperelliptic fibration of
genus 3
$p_2\colon S\to\pp^1$.
\end{itemize}

The singular fibres of $Y\to\pp^1$ consist of a $-1-$curve counted twice (the inverse image
of $P_i'$) and of two disjoint $-2-$curves, contained in the branch locus of $\tilde{S}\to
Y$.  So each singular fibre of $\tilde{S}\to\pp^1$ contains a pair of disjoint $-1-$curves,
 one of which  is contracted by the map $h\colon \tilde{S}\to\F_e$ while the other one is
mapped to a ruling of $\F_e$.  If we denote by $B'$ the branch locus of $h$, then
$B'=D+\Phi_1+\dots +\Phi_6$, where
$\Phi_1\dots
\Phi_6$ are distinct rulings of $\F_e$ and $D$ is effective. In the proof of Theorem 4.4 of
\cite{nodi} it is shown that
 $D$ is linearly equivalent to
$8C_{\infty}+(12+4e)\Phi$. For $i=1\dots 6$ $D$ has a $(4,4)$ point at $P_i\in \Phi_i$ with
the same tangent as $\Phi_i$ and $P_1\dots P_6$ are the only singularities of $B'$.

Notice that the strict transform of $B'$ on $Y$ is  smooth by construction. The
curve $D\subset \F_e$ is the image of the strict transform in $\tilde{S}$ of the divisorial
part of the fixed locus of
$\si$ on
$S$, hence by Theorem \ref{classifica} it has two irreducible components if $(S,\si)$ is of
type Ia and it is irreducible otherwise (cf. the description of the different types of
double planes in \S \ref{sezexamples}).

The map $Y\to\F_e$ is not determined uniquely: two such maps $Y\to \F_e$ and $Y\to\F_{e'}$
are related by  elementary transformations centered at some the points $P_i$ (an
elementary transformation of $\F_e$ centered at $P$ of consists in blowing up $P$ and
then blowing down the strict transform of the ruling through $P$).
If $e=0$, then  we can perform an elementary transformation at, say, $P_1$ and get $e=1$.

The fact that the points $P_i$ are the only singularities of $B'$ and
that they have  ``vertical'' tangent implies that no section of $\F_e$ is contained in
$D$. Thus we have
$0\le C_{\infty}D=12-4e$, namely $e\le 3$ and  at least $3+e$ of the $P_i$  are not in
$C_{\infty}$. Thus, if $e>1$,  one can  perform  an elementary transformation at such a
point  
$P_i$,  replace $e$ by $e-1$ and finally  arrange 
 that $e=1$. 

  If $e=1$,  then the required plane model can be obtained  by composing
$\tilde{S}\to\F_1$ with the birational morphism $\F_1\to \pp^2$ that contracts the
exceptional curve of $\F_1$ and considering the Stein factorization $\tilde{S}\to
X\to\pp^2$. We denote by  $P$  the image of the exceptional curve of $\F_1$,   by $R_i$
the image of
$P_i$,  by $L_i$ the image of $\Phi_i$ and by $C_{16}$ the image of $D$, which is a curve
of degree 16. Since
$C_{\infty}D=8$,
$P$ is a point of multiplicity 8  of $C_{16}$ and  at most two of the
$R_i$ are infinitely near to $P$. The branch
locus of
$X\to\pp^2$ is the image of
$B'$ and it is easy to check that it has the stated properties. In particular, if  $(S,
\si)$ is of type Ia, then the two components of $D$, being numerically equivalent to one
another,  are both mapped to curves of
degree 8.

It is well known  that the
singularities of a double cover of a smooth surface can be solved by repeatedly  blowing
up the base of the cover at the singularities of the branch curve  and taking base change
and normalization, and there  are formulas for the numerical invariants of such a resolution
(cf. for instance
\cite{bpv}, Ch. III, \S 7).
Applying this method to $X\to\pp^2$ one sees that a resolution of $X$ has nonzero geometric
genus iff there exists a curve $C_8\subset \pp^2$ of degree 8 with multiplicity $\ge 6$ at $P$
and with a double point at $R_i$ such that 
$L_i$ is tangent to $C_8$  at $R_i$ for $i=1\dots 6$. The intersection multiplicity of such
a
$C_8$ with a line
$L_i$ is $\ge 9$, hence  $C_8$ must be  equal to $L_1+\dots L_6+C_2$, where
$C_2$ is a conic containing  $R_1\dots R_6$. This proves iv).

A standard computation using again the 
formulas for double covers of \cite{bpv} shows that if $B$ satisfies properties i)--iv)
then the double cover $X\to\pp^2$ is plane model of a minimal  double plane of
general type  $(S,\si)$ with
$p_g=0$ and $K^2=8$ (beware, the resolution of $X$ as a double cover is not minimal!). In
addition,  the pencil of lines through
$P$ pulls  back on $S$  to a hyperelliptic pencil of genus 3, hence $(S,\si)$ is of type I
by Theorem \ref{struttura} and Theorem \ref{classifica}.
\end{proof}
The next theorem describes a plane model for double planes of type II which is quite similar
to the Du Val model for type I.

\begin{thm}\label{planeII}
Let $(S,\si)$ be a minimal double plane of general type with
$K^2=8$ and
$p_g=0$.

If $(S,\si)$ is of type II, then there
 is a plane
model
$X\to\pp^2$ with branch curve $B$ such that:
\begin{itemize}
\item[i)] $B=C_{21}+L_1+\dots + L_5$, where $C_{21}$ has degree 21 and $L_1\dots L_5$ are
distinct lines;
 
\item[ii)] the singularities of
$C_{21}$ are:
\begin{itemize}
\item[--] a  singular  point $P$  of multiplicity 9 such that   $P\in L_1\dots L_5$;

\item[--]  5 points $R_i$ of type $(6,6)$ such that $R_i\in L_i$ and $L_i$ is tangent to
$C_{21}$ at $R_i$,
$i=1\dots 5$;
\end{itemize}

\item[iii)] the singularities of $B$ are
solved by blowing up $P$, $R_1\dots R_5$ and the points $S_i$ infinitely near to $R_i$ in the
direction of $L_i$, $i=1\dots 5$;

\item[iv)] there is no curve $C_5$ of degree 5 containing  $P$ and having a double point at
$R_i$ such that $L_i$ is tangent to $C_5$ at $R_i$ for $i=1\dots 5$.
\end{itemize}

One of the $R_i$ can be infinitely near to $P$. The curve $C_{21}$  is irreducible.

Conversely, given a plane  curve $B$ as in i)--iv), the double cover
$X\to\pp^2$ branched on $B$ is a plane model of a minimal double plane of general type with
$K^2=8$ and $p_g=0$, of type II.
\end{thm}
\begin{proof} The proof is very similar to the proof of Theorem \ref{planeI}, hence it is left
to the reader.
\end{proof}

\section{Moduli}\label{sezmoduli}
In this  section we establish some properties 
of the subset   of the moduli space of surfaces of general type with $p_g=0$ and $K^2=8$
consisting of the surfaces admitting a rational involution. Most of the  arguments that
we use here are  standard in deformation theory, hence  some of proofs are not
very detailed.

We introduce some notation. We denote by $\MM$ be  the moduli space of surfaces of general type with
$p_g=0$ and $K^2=8$ and by $\DD\subset
\MM$  the set of isomorphism classes of surfaces  admitting a rational
involution. We let 
  $\DD_{\rm Ia},\DD_{\rm Ib},\DD_{\rm Ic}, \DD_{\rm Id}, D_{\rm II}\subset
\DD$  be the subsets corresponding to double planes $(S,\si)$ of types Ia,\dots Id, II,
respectively, and we set  $\DD_I:= \DD_{\rm Ia}\cup
\DD_{\rm Ib}\cup
\DD_{\rm Ic}\cup
\DD_{\rm Id}$. 
 Our results are
summarized in the following:

\begin{thm}\label{moduli} Notation and assumptions as above.

Then:
\begin{itemize}
\item[i)]  $\DD$  is a disjoint union:  $\DD=\DD_{\rm Ia}\sqcup
\DD_{\rm Ib}\sqcup
\DD_{\rm Ic}\sqcup
\DD_{\rm Id}\sqcup D_{\rm II}$;
\item [ii)] $D_{\rm I}$ and $D_{\rm II}$ are closed in $\MM$;
\item [iii)]   $\DD_{\rm Ia},\DD_{\rm Ib},\DD_{\rm Ic}, \DD_{\rm
Id}, \DD_{\rm II}$ are normal open and closed subsets of  $\MM$ of the following dimensions:
$$\dim \DD_{\rm Ia}=5,\,\dim \DD_{\rm Ib}=4,\,\dim \DD_{\rm Ic}=3,\,\dim  \DD_{\rm Id}=3,
\, \dim\DD_{\rm II}=2;$$
\item[iv)] $\DD_{\rm Ia},\DD_{\rm Ib},\DD_{\rm Ic}, \DD_{\rm
Id}$ are irreducible.
\end{itemize}
\end{thm} 
Theorem \ref{moduli} describes in particular  the set of surfaces with
non birational bicanonical map. This set is always closed, but in this case it turns out
to be also open. An analogous phenomenon occurs for surfaces with $p_g=0$, $K^2=6$ and
bicanonical map of degree 4 (cf. \cite{mp2}).
\begin{cor}\label{modulibica} The set of surfaces of $\MM$ with non birational
bicanonical map is the union of 4 irreducible connected components of $\MM$ of
respective dimensions 5, 4, 3 and 3.
\end{cor}
\begin{proof}
Follows from Theorem \ref{moduli},  since by Corollary \ref{bica} and  Theorem
\ref{classifica}
 the set of surfaces with non birational bicanonical map is $\DD_{\rm I}$.
\end{proof}
\begin{proof}[Proof of Theorem \ref{moduli}] Statement i) is the content of Corollary
\ref{unique}, ii).  The remaining properties will be a consequence    of the Lemmas
that follow.
\end{proof}

\begin{lem}
The sets $\DD_{\rm I}$ and  $\DD_{\rm II}$ are closed in $\MM$.
\end{lem}
\begin{proof}
By Theorem 4.4 of \cite{nodi} and Theorem \ref{classifica}, $\DD_{\rm I}$ is the set of
surfaces $S$ admitting an involution with 12 isolated fixed points and $\DD_{\rm II}$ is
the set of surfaces admitting an involution with 10 isolated fixed points.  

We remark 
 that by Miyaoka's formula (\cite{miyaoka}, \S 2) a minimal  surface $S$ with $p_g=0$
and $K^2=8$ contains no $-2-$curve, hence it coincides with its canonical model.
Therefore, given a smooth family $\psi
\colon 
\XX\to
\Delta$   
 of minimal surfaces of general type with $K^2=8$ and $p_g=0$, where $\Delta\subset\C$ is
the unit disk, 
 we  have to show that if the class of
$X_t:=\psi\inv(t)$ is in $\DD_{\rm I}$ (resp. $\DD_{\rm II}$) for
each $t\in \Delta\setminus\{0\}$, then also $[X_0]\in \DD_{\rm I}$ (resp. $\DD_{\rm
II}$).

Since by Corollary \ref{unique} the surface $X_t$ admits exactly   one rational involution
$\si_t$, there is a birational map    $\si\colon \XX\to\XX$ that
restricts to 
$\si_t$ on $X_t$ for $t\ne 0$. By Corollary 4.5 of \cite{ritabarbara} $\si$ is actually
biregular and we  
denote by
$\si_0$ the restriction of $\si$ to $X_0$.  Thus  we have to show that the number of
isolated fixed points of $\si_t$ is the same for all $t\in \Delta$.

Clearly one has $\psi\circ\si=\psi$.
Let $P\in \XX$ be a fixed point of $\si$. By Cartan's Lemma and by the fact that $\si$
preserves the fibres of $\psi$, there exist local analytic coordinates near $P$
such that $\psi$ is given by $(x,y,t)\mapsto t$ and $\si$ acts by $(x, y,t)\mapsto
(-x,\epsi y, t)$, where $\epsi=1$ or $\epsi=-1$. Hence the fixed locus of $\si$ is the
disjoint union of two smooth closed sets $\Ga$ and $\Ga'$ of dimensions respectively 1
and 2 and the restriction of $\fie$ to $\Ga$ and $\Ga'$ is a smooth map. So the
$0-$dimensional part of the fixed locus of $\si_t$ on $X_t$ is equal to $X_t\cap \Ga$
for every
$t\in\Delta$. By the above remarks the map  $\Ga\to\Delta$ is proper and \'etale, hence the
cardinality of $X_t\cap \Ga$ is constant.
\end{proof}
The following auxiliary result is likely to be well known, but we give a proof for lack
of a suitable reference.
\begin{lem}\label{omega2}Let $X$ be a smooth variety, let   $G$ be a finite group acting
faithfully on
$X$ in such a way that $Y:=X/G$ is smooth. Let $\pi\colon X\to Y$  be the quotient
map and let $D$  be the reduced branch divisor of $\pi$. Then:
$$(\pi_*\omega_X^2)^G=\omega_Y^2(D).$$
\end{lem}
\begin{proof}
The map $\pi$ is flat and finite. Given a $G-$linearized  line bundle $L$ on $X$, $\pi_*L$
is locally free of   rank equal to $|G|$ and $G$ acts on $\pi_*L$ via the regular
representation. One can define a trace map $tr\colon \pi_*L\to (\pi_*L)^G$\, by
$tr(s):=\frac{1}{|G|}\sum _{g\in G}g^*s$. Clearly $tr$ splits the inclusion
$(\pi_*L)^G\hookrightarrow
\pi_*L$ and thus $(\pi_*L)^G$ is a line bundle. 
If we consider $L=\omega_X^2$ then there is a natural inclusion
$\omega_Y^2\to(\pi_*\omega_X^2)^G$ which is an isomorphism on $Y_0:=Y\setminus D$. So
sections of $(\pi_*\omega_X^2)^G$ are rational sections of $\omega_Y^2$ whose divisor
of  poles is  supported on $D$ and they are characterized by the fact that they pull
back to regular sections of
$\omega_X^2$. The statement follows from this remark by an easy local computation.
\end{proof}
\begin{lem} \label{open}The sets $\DD_{\rm Ia},\DD_{\rm Ib},\DD_{\rm Ic}, \DD_{\rm Id},
\DD_{\rm II}$ are open in $\MM$ and normal.
\end{lem}
\begin{proof}

   For a variety $X$ we denote by
$\Def(X)$ the functor of  deformations of $X$; if there is a group $G$ that acts
on  
$X$ we denote by 
$\Def(X,G)$ the   deformations of $X$ that preserve the $G-$action. It is well known
that the tangent space to $\Def(X,G)$ is $H^1(X,T_X)^G$ and the obstruction space is
$H^2(X,T_X)^G$. 

Consider a point $[S]$ 
 of
$\DD$. By Theorem \ref{classifica},  $S$ is a quotient  
$(F\times C)/G$, where $G$,    $F$ and  $C$ are as explained in \S \ref{sezexamples}. 
 As usual we denote by $G_0$ the subgroup of $\Aut(F)$ generated by $G$
and by  the hyperelliptic involution $\tau$ of $F$.
Consider the map of functors 
$\eta\colon
\Def(F,G_0)\times
\Def(C,G)\to
\Def(S)$ defined by taking fibre products and dividing by the ``diagonal'' $G-$action.
  Clearly,  given an object
$\xi\in 
\Def(F,G_0)\times
\Def(C,G)\to
\Def(S)$,    $\eta(\xi)$ is a family of surfaces whose smooth fibres  have a double plane
structure of the same type as $S$. Since  $F$ and $C$ have dimension 1, 
$\Def(F,G_0)$ and $\Def(C,G)$ are unobstructed. The tangent space to $\Def(S)$ is
$H^1(S,T_S)$. Since the quotient map $F\times C\to S$ is \'etale, $H^1(S,T_S)$ is naturally
isomorphic to
$H^1(F\times C, T_{F\times C})^G$, which is in turn isomorphic to $H^1(F, T_F)^G\oplus
H^1(C,T_C)^G$ by the K\"unneth formula and by the fact that $G$ acts on $F$ and on $C$
separately.  The map on tangent spaces induced by
$\eta$ is the inclusion $H^1(F, T_F)^{G_0}\oplus H^1(C,T_C)^G\to H^1(F, T_F)^G\oplus
H^1(C,T_C)^G$. The claim follows if we show that this map is an isomorphism, i.e. that
$H^1(F, T_F)^{G_0}=H^1(F,T_F)^G$. For types Ia, Ib, Id one has $G=G_0$ and there is nothing
to prove. 
If $S$ is of type Ic or II we  show that $H^1(F,T_F)^G=0$. By Serre duality,
 the dual of  $H^1(F,T_F)$ is $H^0(F,2K_F)$ and
thus it is equivalent to show  $H^0(F,2K_F)^G=0$. In turn this follows by applying Lemma
\ref{omega2} to the map $F\to F/G$, which has three branch points for both types (see 
\S
\ref{sezexamples}).

The above computations also show that the Kuranishi family of $S$ is smooth, hence the
moduli space, which is locally a finite quotient of the Kuranishi family,  is normal.
\end{proof}

\begin{lem}One has:
$$\dim \DD_{\rm Ia}=5,\,\dim \DD_{\rm Ib}=4,\,\dim \DD_{\rm Ic}=3,\,\dim  \DD_{\rm Id}=3,
\,\dim\DD_{\rm II}=2.$$
\end{lem}
\begin{proof}
By the proof of Lemma \ref{open}, if $S=(F\times C)/G$\, is in $\DD$ then the dimension
of
$\DD$ (and of $\MM$) at the point corresponding to $S$ is equal to  the dimension of
$H^1(F,T_F)^G\oplus H^1(C,T_C)^G$.    This
dimension can be computed as in the  proof of Lemma \ref{open}, using Serre duality and
Lemma \ref{omega2}.
\end{proof}

\begin{lem} The sets $\DD_{\rm Ia}, \DD_{\rm Ib}, \DD_{\rm Ic}, \DD_{\rm Id}$ are
irreducible.
\end{lem}
\begin{proof} The proof of Lemma \ref{open} shows that, in order to prove that 
$\DD_{\rm Ia}$ (resp. $\DD_{\rm Ib}, \DD_{\rm Ic},\DD_{\rm Id}$) is irreducible it is
enough to show the existence of two families of smooth curves $\FF\to B_1$ and $\CC\to
B_2$ with a
$G-$action such that:
\begin{itemize}
\item[a)] 
$G$ maps each fibre of $\FF\to B_1$ and $\CC\to B_2$ to itself and the action on the
fibre is the one required for type Ia (resp. Ib, Ic, Id);
\item[b)] every
$F$, resp.
$C$,  with a
$G-$action of the right  type is isomorphic to a fibre of $\FF\to B_1$, resp. 
$\CC\to B_2$.
\end{itemize}
The family $\FF\to B_1$ can be easily constructed using the equations given in \S 
\ref{sezexamples}. We sketch briefly a construction of $\CC\to B_2$ for the various
types. The reader can easily supply the details. For each type we use the notation
introduced in \S
\ref{sezexamples}.
\begin{itemize}
\item[Ia:] the curve $C$ is determined by the three pairs of points $P_{2i-1}+P_{2i}$,
$i=1,2,3$. The only condition is that the six points are distinct;
\item[Ib:] the curve $C$ depends on the
choice of the elliptic curve $E$, of the point $\eta\in E$ of order $4$ and of a pair of
points of $\pp^1=E/D_4$. The only condition is that the two points are not branch points of
the quotient map $E\to E/D_4$;
\item[Ic:]  the curve $C$ is determined by the choice of a homogeneous polynomial of degree
4 with generic branching;
\item[Id:] given the curve $E$, the  map $p\colon E\to\pp^1$ is determined by the two
double fibres, $2P_1+2P_2$ and $2Q_1+2Q_2$. Clearly $\xi:=P_1+P_2-Q_1-Q_2$ is 
$2-$torsion and it is nonzero, since otherwise all the fibres of $p$ over a branch point
would be double. So
$C$ is determined by the choice of $E$, $P_1+P_2$, $\xi$ and of a reduced divisor
$Q_1+Q_2\in |P_1+P_2+\xi|$. The condition that $2Q_1+2Q_2$ and $2P_1+2P_2$ are the only
double fibres   is clearly open, hence also in this case we
have an irreducible parameter space for
$C$.
\end{itemize}
\end{proof}
\begin{rem}
We have not been able to decide whether $\DD_{\rm II}$ is irreducible or not. Actually it
is  not difficult  to give an algebraic construction of a   family  $\CC$ containing all
the $C$ with a $G-$action of the right type,  but it is not clear that it is irreducible.
\end{rem}

\bigskip
\bigskip

\begin{minipage}{13cm}
\parbox[t]{8cm}{Rita Pardini\\
Dipartimento di Matematica\\
Universit\`a di Pisa\\
Via Buonarroti, 2\\
56127 Pisa, ITALY\\
pardini@dm.unipi.it}
\end{minipage}
\end{document}